\documentclass[11pt,twoside]{article}
\usepackage{amsfonts}
\usepackage{color}
\usepackage{fancyhdr}
\usepackage{titlesec}
\usepackage{cite}
\usepackage{ifthen}
\usepackage{amssymb}
\usepackage{fancyhdr}
\usepackage{titlesec}
\usepackage{pifont}
\usepackage{setspace}
\usepackage{indentfirst}
\usepackage{amsmath,amssymb,amscd,amsthm,mathrsfs}
\input amssym.def

\newboolean{first}
\setboolean{first}{true}
\renewcommand{\headrulewidth}{0pt}

\newfont{\aaa}{cmb10 at 19pt}
\newfont{\bbb}{cmb10 at 11pt}
\newtheorem{lemma}{Lemma}[section]
\newtheorem{theorem}{Theorem}[section]

\newtheorem{rem}{Remark}[section]

\newtheorem{proposition}{Proposition}[section]
\newcommand{\beq}{\begin{equation}}
\newcommand{\eeq}{\end{equation}}
\newcommand{\bey}{\begin{eqnarray}}
\newcommand{\eey}{\end{eqnarray}}
\newcommand{\beyy}{\begin{eqnarray*}}
\newcommand{\eeyy}{\end{eqnarray*}}

\numberwithin{equation}{section}
\makeatletter
\def\@evenfoot{}
\def\@oddfoot{}
\makeatother
\makeatletter
\newcommand{\Rmnum}[1]{\expandafter\@slowromancap\romannumeral #1@}
\makeatother

\begin{document}
\thispagestyle{empty} \thispagestyle{fancy} {
\fancyhead[RO,LE]{\scriptsize \bf 
} \fancyfoot[CE,CO]{}}
\renewcommand{\headrulewidth}{0pt}
\begin{center}
{\bf \LARGE A generalization of a lemma of Boccardo and Orsina and application}

\vspace{2mm}

{\small \textsc{Hongya GAO}$^1$  \qquad   \textsc {Meng GAO}$^1$  \qquad   \textsc {Siyu GAO}$^2$\footnote{Corresponding author, email: siyugao@my.unt.edu.}

{\small 1. College of Mathematics and Information Science, Hebei University, Baoding, 071002, China\\
2. Department of Mathematics, University of North Texas, Denton, Texas 76203, USA}}
\end{center}

\vspace{3mm}

\begin{center}
	\begin{minipage}{130mm}
		
		{\bf \small Abstract}.  {\small We present a generalization of a technical lemma due to Boccardo and Orsina, and then give an application to regularity of minima for integral functionals noncoercive in the energy space.}
					
		
	\end{minipage}
\end{center}

\thispagestyle{fancyplain} \fancyhead{}
\fancyhead[L]{\textit{}\\
} \fancyfoot{}

\section{Introduction.}

In dealing with regularity properties of minima of some integral functionals  noncoercive in the energy space, Boccardo and Orsina proved in \cite {BO} a useful lemma to consider how the regularity of $u$ depends on the summability of the source $f$ in some Marcinkiewicz space. More precisely, Boccardo and Orsina considered integral functionals of the type
\begin{equation}\label{integral functional No.1}
{\cal J}(v) =\int_\Omega a(x,v) j(\nabla v) dx -\int_\Omega fv dx, \ \ v\in W_0^{1,p} (\Omega),
\end{equation}
where $\Omega $ is a bounded open subset of $\mathbb R^n$, $n\ge 2$,
\begin{equation}\label{condition for f}
f\in L^r(\Omega),  \ r\ge (p^*)',
\end{equation}
$a(x,s):\Omega \times \mathbb R \rightarrow \mathbb R$ is a Carath\'eodory function (that is, $a(\cdot, s)$ is measurable on $\Omega$ for every $s$ in $\mathbb R$, and $a(x,\cdot)$ is continuous on $\mathbb R$ for almost every $x$ in $\Omega$) such that for almost every $x\in \Omega $ and every $s\in \mathbb R$,
\begin{equation}\label{condition for a(x) No.1}
a(x,s) =\frac {\beta_1}{ (b(x)+|s|) ^{\alpha p}},
\end{equation}
where $\beta_1>0$, $1<p <n$,
\begin{equation}\label{condition for alpha No. 1}
 0<\alpha < \frac 1 {p'},
\end{equation}
and $b(x)$ is a measurable function on $\Omega $ such that
\begin{equation}\label{condition for b}
 0<\beta _2 \le \beta (x) \le \beta_3 <+\infty, \ \ \mbox { for almost every } x\in \Omega;
\end{equation}
moreover, $j:\mathbb R^n \rightarrow \mathbb R$ is a convex function such that $j(0)=0$ and for every $\xi \in \mathbb R^n$,
\begin{equation}\label{condition for j No. 1}
\beta_4 |\xi| ^p \le j(\xi)\le \beta_5(1+|\xi| ^p)
\end{equation}
with $\beta _4,\beta _5>0$.

Note that, although the functional (\ref{integral functional No.1}) is well-defined for $f$ satisfying (\ref{condition for f}), it is not coercive on the energy space $W_0^{1,p} (\Omega)$: there exists a function $f$, and a sequence $\{u_n\}$ whose norm diverges in $W_0^{1,p} (\Omega)$, such that ${\cal J} (u_n)$ tends to $-\infty$, see Example 3.3 in \cite{BO}. Thus $\cal J$ may not attain its minimum on $W_0^{1,p} (\Omega)$. A good idea to consider a minimum of $\cal J$ is to extend $W_0^{1,p} (\Omega)$ to a larger space $W_0^{1,q} (\Omega)$ with
\begin{equation}\label{definitioin for q}
q=\frac {np (1-\alpha)}{n-\alpha p} <p,
\end{equation}
in the following way:
\begin{equation}\label{definition for I}
{\cal I} (v) =\left \{
\begin{array}{llll}
{\cal J} (v), & \mbox { if  } {\cal J}(v) \mbox { is finite} , \\
+\infty, & \mbox { otherwise} .
\end{array}
\right.
\end{equation}

A function $u\in W_0^{1,q} (\Omega)$ is called a minima of the functional ${\cal I}(v)$ in (\ref{definition for I}) if
\begin{equation}\label{definition for minima}
{\cal I}(u) \le {\cal I} (v) , \ \  \mbox { for all } v\in W_0^{1,q} (\Omega).
\end{equation}

The following result comes from Theorem 2.1 in \cite{BO}:
Suppose that $f$ belongs to $L^r(\Omega)$, with $r\ge [p^* (1-\alpha)]'$. Then $\cal I$ is coercive and weakly lower semicontinuous on $W_0^{1,q} (\Omega)$. By standard results, $\cal I$ has a minimum on $W_0^{1,q} (\Omega)$. See Theorem 1.1 in \cite{BO}.

Let us  recall the Marcinkiewicz space $M^{r}(\Omega)$, $r>0$, which is the set of all measurable functions $f: \Omega\rightarrow \mathbb R$ such that
\begin{equation}\label{definition for M space}
\left| \{x\in\Omega: |f(x)|>t\} \right| \leq \frac{c}{t^{r}},
\end{equation}
for every $t>0$ and for some constant $c>0$. The norm of $f\in M^r (\Omega)$ is defined by
$$
\|f\| _{M^r (\Omega)} ^r =\inf \left \{c>0:  (\ref{definition for M space}) \ \rm { holds} \right\}.
$$
If $\Omega$ has finite measure, then
\begin{equation}\label{LML}
L^{r}(\Omega)\subset M^{r}(\Omega)\subset L^{r-\varepsilon}(\Omega),
\end{equation}
for every $r>1$ and every $0<\varepsilon\leq r-1$. We recall also that, see Proposition 3.13 in \cite{Boccardo-Croce}, if $f\in M^{r}(\Omega)$, $r>1$, then there exists a positive constant $B=B(\|f\| _{M^r(\Omega)},r)$ such that, for every measurable set $E\subset \Omega$,
\begin{equation}\label{Holder}
\int_{E}|f|dx\leq B |E|^{1-\frac{1}{r}}.
\end{equation}

In \cite{BO}, Boccardo and Orsina obtained some regularity results for $u$ and $\nabla u$ in terms of
\begin{equation}\label{condition for f No.2}
f\in M^r(\Omega), \ \ [p^*(1-\alpha)]' <r<\frac n p,
\end{equation}
see Theorem 6.3 in \cite{BO}.

\begin{proposition}\label{Boccardo and Orsina}
Let $f$ be in $M^r(\Omega)$, with $[p^*(1-\alpha)]' <r<\frac n p$. Then any minima $u$ of $\cal I $ belongs to $M^s (\Omega)$ with
$$
s =\frac {nr (p(1-\alpha) -1 )}{n-rp}.
$$
Moreover:

{\rm a)} if
$$
\left(\frac {p^*}{1+\alpha p}\right) ' <r<\frac n p,
$$
then $u$ belongs to $W_0^{1,p} (\Omega)$;

{\rm b)} if
$$
[p^*(1-\alpha)] ' < r  \le \left(\frac {p^*}{1+\alpha p}\right) ',
$$
then $|\nabla u| $ belongs to $M^\rho (\Omega)$ with
$$
\rho = \frac {nr [p(1-\alpha ) -1]}{n-r (1+\alpha p )}.
$$
 \end{proposition}

The main tool in proving Proposition \ref{Boccardo and Orsina} is a technical lemma as follows, see Lemma 6.1 in \cite{BO}.

\begin{lemma}\label{BO}
Let $\psi :[0,+\infty) \rightarrow [0,+\infty) $ be a non increasing function, and suppose that
\begin{equation}\label{the basic inequality}
\psi (h) \le c \frac {k^A \psi (k) ^B +\psi (k) ^C}{ (h-k)^D}, \ \ \forall h>k\ge 0,
\end{equation}
where $c$ is a positive constant, and $A,B,C$ and $D$ are such that
\begin{equation}\label{conditions for ABCD}
A<D, \ C<B<1, \ \frac {D-A}{1-B} =\frac {D}{1-C}.
\end{equation}
Then there exists $\bar k \ge 0$, and a positive constant $\bar c $ such that
\begin{equation}\label{conditions for ABCD}
\psi (k) \le \bar c k ^{-\frac {D-A}{1-B}} =\bar c k ^{-\frac {D}{1-C}}, \ \ \forall k\ge \bar k.
\end{equation}
\end{lemma}

We remark that this lemma is very similar to the classical Stampacchia Lemma (see \cite{Stampacchia}, Lemma 4.1) used till now repeatedly by many mathematicians in dealing with regularity issues of solutions of elliptic PDEs as well as minima of variational integrals.

Let us come back to Proposition \ref{Boccardo and Orsina}. For $v\in W_0^{1,q} (\Omega)$ and $f\in M^r(\Omega)$, $r > [p^*(1-\alpha)]' =( q^*)'$, the second integral  $\int_\Omega fv dx $ is well-defined thanks to Sobolev embedding. Two natural questions now arise as:  do we have any regularity properties for minima of $\cal I$ when $r=\frac n p $ and $r>\frac n p$ ? In order to answer these two questions, we need to generalize Lemma \ref{BO}.

\section{A generalization of Lemma \ref{BO}.}

We now prove the following

\begin{lemma}\label{generalization of Stampacchia Lemma}
Suppose $c_1, A, B, C, D$ are positive constants, $A<D$, $k_0 \ge 0$. Let $\psi :[k_0,+\infty) \rightarrow [0,+\infty) $ be a non increasing function, and suppose that
	\begin{equation}\label{the basic inequality  2}
		\psi (h) \le c_1 \frac {h^A \psi (k) ^B +\psi (k) ^C}{ (h-k)^D},
	\end{equation}
for every $h, k$ with $h>k>k_0 $. It results that:

{\bf i)} if $C<B<1$ and $\frac {D-A}{1-B} =\frac {D}{1-C}$, then there exists a positive constant $\bar c_1$ such that  for any $k\ge k_0$,
\begin{equation} \label{C<B<1}
	\psi (k) \le \bar c_1 k ^{-\frac {D-A}{1-B}} = \bar c_1 k ^{-\frac {D}{1-C}};
\end{equation}

{\bf ii)} if $B=C=1$, then for any $k\ge k_0$ we have
\begin{equation}
	\label{B=C=1}
	\psi (k) \le \psi(k_0) e^{1-\left( \frac{k-k_0}{\tau}\right) ^{\frac{D-A}{D}}},
\end{equation}	
where
\begin{equation}
	\label{B=C=1,(ii)}
\tau =\max \left\lbrace k_0+1, \left(\dfrac{2 c_1 e2^{\frac{(2D-A)A}{D-A}}(D-A)^D}{D^D}\right)^{\frac{1}{D-A}}  \right\rbrace ;
\end{equation}	

{\bf iii)} if $B>C>1$, then
\begin{equation}\label{B>C>1}
\psi(2L) = 0,
\end{equation}
\noindent
where
\begin{equation}\label{L}
\begin{array}{llll}
L=&\displaystyle \max \Big\{1, 2k_0, \left(c_12^{1+D} (1+\psi (k_0)) ^B \right)^{\frac 1 {D-A}}, \\
&\displaystyle  \left( c_1 ^{\frac{C}{C-1}}(1+ \psi(k_0))^{B}2^{D+1 +\frac{A+D+1 }{C-1}+\frac{D}{(C-1)^2}} \right) ^{\frac {C-1}{ (D-A)C}}\Big\} .
\end{array}
\end{equation}
\end{lemma}

We mention that the difference between Lemma \ref{generalization of Stampacchia Lemma} i) and Lemma \ref{BO} is that a real number $k_0$ replaced $0$ and $h^A $ replaced $k^A$. We mention also that, this section borrows some ideas from the papers \cite{BO,S,KV,GLW,GDHR,GHR,GZM}.

\vspace{2mm}

In the proof of Lemma \ref{generalization of Stampacchia Lemma} iii) we shall use the following lemma, which comes from Lemma 7.1 in\cite{E}.

\begin{lemma}\label{E1}
Let $\beta, M, \bar C, x_i$ be such that $\beta >1$, $\bar C>0$, $M>1$, $x_i\geq 0$ and
\begin{equation}
	\label{equatipn-E}
	x_{i+1}\leq \bar C M^{i}x_{i}^{\beta},~~~i=0,1,2,\cdots .
\end{equation}
If
$$x_0\leq \bar C^{-\frac{1}{\beta-1}}M^{-\frac{1}{(\beta-1)^2}},$$
then, we have
$$x_i\leq M^{-\frac{i}{\alpha}}x_0,~~~i=0,1,2,\cdots,$$
so that
$$\lim _{i\rightarrow +\infty} x_i =0.$$
\end{lemma}

\noindent {\bf Proof of Lemma \ref{generalization of Stampacchia Lemma}}.

{\bf i)} Define
\begin{equation}\label{definition for lambda}
\lambda=\frac{D-A}{1-B} =\frac {D}{1-C}
\end{equation}
and
$$
\rho(h)=h^{\lambda}\psi(h).
$$
Then (\ref{the basic inequality 2}) implies
$$\rho(h)\leq c_1\frac{h^{\lambda}\left( h^{A}\psi(k)^{B}+\psi(k)^{C}\right) }{(h-k)^{D}}, \ \ \forall h>k\ge k_0.
$$
Choosing $h=2k$ in the above inequality and we have, for all positive $k\ge k_0$,
\begin{equation*}
\begin{array}{llll}
	\displaystyle \rho(2k)& \leq \displaystyle  c_1\frac{(2k)^{\lambda}\left( (2k)^{A}\psi(k)^{B}+\psi(k)^{C}\right) }{k^{D}}\\
	& \leq \displaystyle c_1\frac{2^{\lambda+A}\left( (k^{\lambda+A}\psi(k)^{B}+k^\lambda\psi(k)^{C}\right) }{k^{D}}\\
	& =\displaystyle c_1\frac{2^{\lambda+A}\left( (k^{\lambda+A-\lambda B}\rho(k)^{B}+k^{\lambda-\lambda C}\rho(k)^{C}\right) }{k^{D}}.\\
\end{array}
\end{equation*}
From the definition of $\lambda $ in (\ref{definition for lambda}), one has $\lambda+A-\lambda B=D$ and $\lambda-\lambda C=D$, then the above inequality becomes
\begin{equation}\label{2.8}
\rho(2k)\leq c_12^{\lambda+A}\left( \rho(k)^B+\rho(k)^C\right), \  \mbox { for all positive } k \ge k_0.
\end{equation}

We claim that, for every integer $n\ge 0$,
 \begin{equation}\label{2^n k0}
 	\rho(2^{n} k_0 )\leq c_1^{\frac 1 {1-B}} 2^ {\frac{\lambda+A+1}{1-B}}\left( 1+\rho( k_0)\right) ^{B^n}.
 \end{equation}
For $n=0$, (\ref{2^n k0}) is obvious (it is no loss of generality to assume $c_1\ge 1$ so that $c_1^{\frac 1 {1-B}} 2^ {\frac{\lambda+A+1}{1-B}} \ge 1$, which ensures (\ref{2^n k0}) in case of $n=0$). We now suppose, for some $n \in \mathbb N$, (\ref{2^n k0}) holds true, and we proceed by induction on $n$. Since $C<B<1$ , then we have,  using (\ref{2.8}) with $k= k_0  $ that
 $$
 \begin{array}{llll}
 \displaystyle \rho (2  k_0  ) &\displaystyle  \leq c_12^{\lambda+A}\left( \rho(  k_0   )^B+\rho(  k_0   )^C\right)  \\
 & \displaystyle  \leq  c_12^{\lambda+A+1}(1+\rho(  k_0   ))^B   \\
&\displaystyle  \leq  c^{\frac 1 {1-B}}_12^{\frac{\lambda+A+1}{1-B}}(1+\rho(  k_0   ))^B.
 \end{array}
 $$
 By (\ref{2.8}) again, (\ref{2^n k0}) and the above inequality that
 \begin{equation*}
 	\begin{array}{llll}
 		\displaystyle\rho(2^{n+1}  k_0   )&\leq\displaystyle c_12^{\lambda+A}\left( \rho(2^{n}  k_0   )^B+\rho(2^{n}  k_0   )^C\right) \\
 		&\leq \displaystyle c_12^{\lambda+A}\left[ \left(c ^{\frac 1 {1-B}}_1 2^{\frac{\lambda+A+1}{1-B}}(1+\rho(  k_0   ))^{B^n}\right) ^{B}+\left(c_1^{\frac 1 {1-B}} 2^{\frac{\lambda+A+1}{1-B}}(1+\rho(  k_0   ))^{B^n}\right)^{C}\right] \\
 		&\leq \displaystyle c_1 ^{1+\frac {B}{1-B}}2^{\lambda+A}\left[ 2^{\frac{(\lambda+A+1)B}{1-B}}(1+\rho(  k_0   ))^{B^{n+1}}+2^{\frac{(\lambda+A+1)B}{1-B}}(1+\rho(  k_0   ))^{B^{n+1}}\right] \\
 		&\le c ^{\frac 1 {1-B}}_12^{(\lambda+A+1)(1+\frac{B}{1-B})}(1+\rho(  k_0   ))^{B^{n+1}}\\
 		&=c ^{\frac 1 {1-B}}_12^{\frac{\lambda+A+1}{1-B}}(1+\rho(  k_0   ))^{B^{n+1}},
 	\end{array}
 \end{equation*}
 which is (\ref{2^n k0}) for $n+1$. Thus, (\ref{2^n k0}) holds for every $n\ge 0$.

 Since $B<1$, then for $n\ge 0$, $(1+\rho ( k_0    ))^{B^n}\le (1+\rho ( k_0   ))^B$, from (\ref{2^n k0}) it follows that
\begin{equation}\label{definition for M}
 \rho(2^{n} k_0   )\leq c ^{\frac 1 {1-B}}_12^{\frac{\lambda+A+1}{1-B}}(1+\rho( k_0   ))^B:=M, \ \ n\ge 0,
\end{equation}
 and so, recalling the definition of $\rho$,
 \begin{equation}\label{2^n k0 1}
 	\psi(2^{n}  k_0   )\leq \frac{M}{(2^{n}  k_0   )^{\lambda}}.
 \end{equation}
 For any positive $k\ge  k_0  $, there exist $k'\in[ k_0   ,2  k_0   )$ and $n\ge 0$ such that $k=2^{n}k'$, and so $2^{n}  k_0   \leq k < 2^{n+1}  k_0   $. Since $\psi$ is non increasing, we thus have, by (\ref{2^n k0 1}),
 $$
 \psi(k)\leq \psi(2^{n} k_0 )\leq \frac{M}{(2^{n}  k_0   )^{\lambda}}=  \frac{2^{\lambda}M}{(2^{n+1}  k_0   )^{\lambda}}\le \frac{2^{\lambda}M}{k^{\lambda}},
 $$
 so that Claim i) is proved with $\bar c_1 =2^\lambda M $ with $\lambda $ and $M $ be as in (\ref{definition for lambda}) and (\ref{definition for M}) respectively.

{\bf ii)} Let $B=C=1$ and $\tau$ be as in (\ref{B=C=1,(ii)}). For $s=0,1,2,\cdots $,  we let
$$
k_s=k_0+\tau s^{\frac{D}{D-A}},
$$
then $\{k_s\}$ is an increasing sequence and
$$
k_{s+1}-k_{s}=\tau \left[ (s+1)^{\frac{D}{D-A}}-s^{\frac{D}{D-A}} \right].
$$
We use Taylor's formula  to get
\begin{equation}\label{Taylor}
k_{s+1}-k_{s}=\tau\left[  \frac{D}{D-A}s^{\frac{A}{D-A}}+\frac{AD}{(D-A)^2}\xi^{\frac{2A-D}{D-A}}\right]  \geq\frac{\tau D}{D-A}s^{\frac{A}{D-A}},
\end{equation}
where $\xi$ lies in the open interval $(s,s+1)$. In (\ref{the basic inequality 2}) we take $B=C=1$, $k=k_s$ and $h=k_{s+1}$, we use (\ref{Taylor}) and we get, for $s \in \mathbb N^+ =\{1,2,\cdots\}$,
\begin{equation}\label{k_s+1}
\begin{array}{llll}
\psi(k_{s+1})&\leq c_1\displaystyle \dfrac{\left[ k_0+\tau (s+1)^{\frac{D}{D-A}}\right] ^{A}\psi(k_s)+\psi(k_s)}{\left( \frac{\tau D}{D-A}\right) ^{D}s^{\frac{AD}{D-A}}}\\
&\leq 2 c_1 \dfrac{\left[ k_0+1+\tau(2s)^{\frac{D}{D-A}}\right] ^{A}\psi(k_s)}{\left( \frac{\tau D}{D-A}\right) ^{D}s^{\frac{AD}{D-A}}}.
\end{array}
\end{equation}
(\ref{B=C=1,(ii)}) ensures, for $s\in \mathbb N^+$,
$$
k_0+1\leq\tau<\tau(2s)^{\frac{D}{D-A}}~~\mbox { and } ~~\dfrac{2 c_1\left( 2^{1+\frac{D}{D-A}}\tau\right) ^{A}}{\left( \frac{\tau D}{D-A}\right) ^{D}}\leq\frac{1}{e}.
$$
From (\ref{k_s+1}) and the above inequalities, one has
\begin{equation*}
\psi(k_{s+1}) \leq \dfrac{ 2 c_1\left( 2\tau(2s)^{\frac{D}{D-A}}\right) ^{A}}{\left( \frac{\tau D}{D-A}\right) ^{D}s^{\frac{AD}{D-A}}}\psi(k_s)
		=\dfrac{2 c_1\left( 2^{1+\frac{D}{D-A}}\tau\right) ^{A}}{\left( \frac{\tau D}{D-A}\right) ^{D}}\psi(k_s)
		\leq\frac{1}{e}\psi(k_s).
\end{equation*}
By recursion,
$$
\psi(k_s)\leq\frac{1}{e^s}\psi(k_0),~~~\forall s\in\mathbb N^+.
$$
The above inequality holds true for $n=0$ as well.
For any $k\geq k_0$, there exists $s\in \{0,1,2,\cdots, \}$ such that
$$
k_0+\tau s ^{\frac{D}{D-A}}\leq k<k_0+\tau (s+1)^{\frac{D}{D-A}}.
$$
Thus, considering $\psi(k)$ is nonincreasing, one has
$$
\psi(k)\leq \psi\left( k_0+\tau s^{\frac{D}{D-A}}\right) =\psi(k_{s})\leq e^{-s}\psi(k_0)\leq \psi(k_0)e^{1-\left( \frac{k-k_0}{\tau}\right) ^{\frac{D-A}{D}}},
$$
as desired.

{\bf iii)} For $B>C>1$ we fix $L\ge \max \{1,2 k_0 \}$ (which implies $L-k_0 \ge \frac L 2$) such that $\psi(L)\leq 1$. This can always be done since one can choose in (\ref{the basic inequality  2}) $h=L$, $k=k_0$, using the fact $A<D$, one has
$$
\psi (L) \le c_1 \frac {L^A \psi (k_0)^B +\psi (k_0)^C}{(L-k_0 ) ^D} \le c_1 \frac {2L^A (1+\psi (k_0)) ^B}{ (L-k_0) ^D} \le c_1 \frac {2 ^{1+D} (1+\psi (k_0)) ^B}{ L^{D-A}},
$$
thus $\psi (L)\le 1 $ would be satisfied if
\begin{equation}\label{psi L ge 1}
L^{D-A} \ge c_1 2 ^{1+D} (1+\psi (k_0)) ^B.
\end{equation}

We choose levels
$$
k_i=2L(1-2^{-i-1}),~~~i=0,1,2,\cdots.
$$
It is obvious that
$$
k_0=L\leq k_i<2L,
$$
$\{k_i\}$ be an increasing sequence, and
$$
\lim _{i\rightarrow +\infty} k_i =2L.
$$
We choose in (\ref{the basic inequality 2})
$$
k=k_i,~  h=k_{i+1},
$$
let
$$
x_i=\psi(k_i),\  x_{i+1}=\psi(k_{i+1}),
$$
and notice that
$$
h-k=k_{i+1}-k_i=L2^{-i-1}, \ \ x_i=\psi (k_i)\le \psi (L) \le 1,
$$
we have, for $i=0,1,2,\cdots$,
\begin{equation*}
\begin{array}{llll}
x_{i+1}&\leq \displaystyle c_1 \dfrac{\left[ 2L(1-2^{-i-2})\right] ^{A}x_i^B+x_i^C}{(L2^{-i-1})^{D}}\\
&\leq \displaystyle 2c_1\frac{\left[ 2L(1-2^{-i-2})\right] ^{A}}{(L2^{-i-1})^{D}} x_i^C\\
&\leq \displaystyle c_1\frac{2^{A+D+1}(2^D)^{i} }{L^{D-A}} x_i^C.
\end{array}
\end{equation*}
Thus (\ref{equatipn-E}) holds true with
$$
\bar C=\frac{c_1 2^{A+D+1 }}{L^{D-A}},  ~M=2^D ~ \mbox { and }  ~ \beta =C>1.
$$
We use Lemma \ref{E1} and we have
\begin{equation}\label{rightarrow 0}
\lim _{i\rightarrow +\infty} x_i =  \lim _{i\rightarrow +\infty} \psi (k_i) =0
\end{equation}
provided that
\begin{equation}\label{condition}
x_0=\psi(k_0)=\psi(L)\leq \left( \frac{c_12^{A+D+1 }}{L^{D-A}}\right) ^{-\frac{1}{C-1}}(2^D)^{-\frac{1}{(C-1)^2}}.
\end{equation}
Note that (\ref{rightarrow 0}) implies
$$
\psi(2L)=0.
$$
Let us check condition (\ref{condition}) and determine the value of $L$. (\ref{condition}) is equivalent to
\begin{equation}\label{equivalent condition}
\psi(L)\leq c_1 ^{-\frac{1}{C-1}}2^{-\frac{A+D+1 }{C-1}-\frac{D}{(C-1)^2}}L^{\frac{D-A}{C-1}}.
\end{equation}
In (\ref{the basic inequality 2}) we take $k=k_0$ and $h=L\geq \max \{1, 2k_{0} \}$ and we have
$$
\begin{array}{llll}
\displaystyle \psi(L) & \displaystyle \leq c_1\dfrac{L^A\psi(k_0)^B +\psi(k_0)^C}{(L-k_0)^D}  \leq \frac{2 c_1L^A (1+\psi(k_0))^B}{(L-k_0)^D} \\
& \displaystyle \leq \frac{2^{D+1}c_1 L^A (1+\psi(k_0))^B }{L^D} =\frac{2^{D+1}c_1 (1+\psi(k_0)) ^B}{L^{D-A}}.
\end{array}
$$
Then (\ref{equivalent condition}) would be satisfied if  $\psi (L) \le 1$,
\begin{equation}\label{condition1}
L\geq \max \{1, 2k_0\} ,
\end{equation}
and
$$
\frac{2^{D+1}c_1 (1+\psi(k_0)) ^B}{L^{D-A}} \leq c_1 ^{-\frac{1}{C-1}}2^{-\frac{A+D+1 }{C-1}-\frac{D}{(C-1)^2}}L^{\frac{D-A}{C-1}}.
$$
The above inequality is equivalent to
\begin{equation}\label{condition2}
c_1 ^{\frac{C}{C-1}}(1+ \psi(k_0))^{B}2^{D+1 +\frac{A+D+1 }{C-1}+\frac{D}{(C-1)^2}}\leq L^{\frac{(D-A)C}{C-1}}.
\end{equation}
(\ref{L}) is a sufficient condition for (\ref{psi L ge 1}), (\ref{condition1}) and (\ref{condition2}).

This ends the proof of Lemma \ref{generalization of Stampacchia Lemma}.   \qed

\vspace{3mm}

We notice that, in the proof of Claim i) of Lemma \ref{generalization of Stampacchia Lemma} we have chosen $h=2k$. Let us take $h=2k$ in (\ref{the basic inequality  2}) and replace $c_1$ by $c_2$, that is,
\begin{equation}\label{the basic inequality 3}
			\psi (2k) \le c_2 \frac {(2k)^A \psi (k) ^B +\psi (k) ^C}{ k^D}, \ \ \forall k\ge k_0.
		\end{equation}
\begin{rem}\label{remark 1} Assume (\ref{the basic inequality  2}). We take $h=2k$ and we get (\ref{the basic inequality 3}) with $c_2=c_1$. That is, (\ref{the basic inequality  2}) implies (\ref{the basic inequality 3}).
\end{rem}

\begin{rem}\label{remark 2}
If we replace (\ref{the basic inequality  2}) by (\ref{the basic inequality 3}), then the Claim i) of Lemma \ref{generalization of Stampacchia Lemma} also holds.
\end{rem}

In the remaining part of this section, we shall consider the relationship between (\ref{the basic inequality 3}) and (\ref{the basic inequality  2}).
We ask the following question:
is (\ref{the basic inequality 3}) weaker than (\ref{the basic inequality  2})? The answer is: it depends. For different cases, we have different answers. We give the following three remarks.


\begin{rem}
For the case $C<B<1$, one has
$$
	(\ref{the basic inequality  2})\Leftrightarrow(\ref{the basic inequality 3}).
$$
\end{rem}
\noindent {\bf Proof. }     ``$\Rightarrow$".  See Remark \ref{remark 1}.

``$\Leftarrow$". Assume (\ref{the basic inequality 3}). Let us consider $h>k\geq k_0$. We split the proof into two cases: $2^{n+1}k\geq h>2^n k$ for some integer $n\geq1$ and $2k\geq h>k$.

{\bf Case $2^{n+1}k\geq h>2^n k$ for some integer $n\geq1$.} Since $\psi$ decreases, we have $\psi(h)\leq \psi(2^n k)=\psi(2(2^{n-1}k))$; we keep in mind that $n\geq1$ so $2^{n-1}k\geq k\geq k_0$ and we can use (\ref{the basic inequality 3}) with $2^{n-1}k$ in place of $k$: we have
$$
\begin{array}{llll}
\displaystyle \psi(2(2^{n-1}k)) & \displaystyle  \leq c_2\dfrac{(2^nk)^A\left[ \psi(2^{n-1}k)\right] ^B+\left[ \psi(2^{n-1}k)\right] ^C}{(2^{n-1}k)^{D}} \\
& \displaystyle  \leq c_2\dfrac{h^A\left[ \psi(2^{n-1}k)\right] ^B+\left[ \psi(2^{n-1}k)\right] ^C}{(2^{n-1}k)^{D}}.
\end{array}
$$
Since $2^{n-1}k\geq k$, we use the monotonicity of $\psi$ to have $\psi(2^{n-1}k)\leq \psi(k)$; then $\left[ \psi(2^{n-1}k)\right] ^{B}\leq \left[ \psi(k)\right] ^{B}$ and $\left[ \psi(2^{n-1}k)\right] ^{C}\leq \left[ \psi(k)\right] ^{C}$. Since $2^{n+1}k\geq h$, we have $(2^{n+1}-1)k\geq h-k$, then
$$
2^{n-1}k=\frac{2^{n+1}k}{4}\geq \frac{(2^{n+1}-1)k}{4}\geq \frac{h-k}{4},
$$
thus
\begin{equation*}
	\begin{array}{llll}
		\psi(h) & \leq\psi(2^n k)=\psi(2(2^{n-1}k)) \\
&\leq c_2\dfrac{h^{A}\left[ \psi(2^{n-1}k)\right] ^{B}+\left[ \psi(2^{n-1}k)\right] ^{C}}{(2^{n-1}k)^{D}}\\
		&\leq 4^D c_2 \displaystyle\frac{ h^A \left[ \psi(k)\right] ^B+ \left[ \psi(k)\right] ^C}{(h-k)^D}.
	\end{array}
\end{equation*}
{\bf Case $2k\geq h>k$.} We use Remark \ref{remark 2} and Claim i) of Lemma \ref{generalization of Stampacchia Lemma} and we get (\ref{C<B<1}).
  We use (\ref{C<B<1}) and the fact that $\psi$ decreases in order to get
\begin{equation*}
\begin{array}{llll}
&  \psi(h) \le \psi(k)=\left[ \psi(k)\right] ^B \left[ \psi(k)\right]^{1-B} \\
&\leq \left[ \psi(k)\right] ^B\left[\bar c_1 \left( \frac 1 k \right) ^{\frac {D-A}{1-B}}\right] ^{1-B}
=\displaystyle \bar c_1^{1-B}k^A \left[ \psi(k)\right] ^B
\left( \frac{1}{k}\right) ^{D}.
\end{array}
\end{equation*}
Since $2k\geq h$ we get $k\ge h-k$ and $\left( \frac{1}{k}\right) ^D \le \left( \frac{1}{h-k}\right) ^D$, then
$$
\begin{array}{llll}
\psi(h) &\displaystyle \leq \frac{c_1 ^{1-B} h^A \left[ \psi(k)\right] ^B }{(h-k)^D}
 &\displaystyle \leq \bar c_1^{1-B}\frac{h^A\left[ \psi(k)\right] ^B+\left[ \psi(k)\right] ^C}{(h-k)^D}.
\end{array}
$$
In both cases we have obtained (\ref{generalization of Stampacchia Lemma}) with $c_1=\max \Big\{ 4^D c_2;  \bar c_1 ^{1-B}\Big \}$.   \qed

\begin{rem} For the case $B=C=1$, one has
	\begin{equation*}\label{R2}
		(\ref{the basic inequality  2})\nLeftarrow(\ref{the basic inequality 3}).
	\end{equation*}
More precisely, the function
\begin{equation}\label{example 11}
	\psi(k)=e^{-( \ln k)^2 },~~~k\in[1,+\infty)
\end{equation}
verifies (\ref{the basic inequality 3}) with $k_0=1$, $B=C=1$, $c_2=\frac 1 {2\ln 2}$, $D=2\ln 2$, and any $0< A< 2\ln2$, but it does not satisfy (\ref{the basic inequality  2}) with $B=C=1$, for any choice of the three constants $D>A>0$ and $c_1>0$.
\end{rem}

\noindent {\bf Proof. }
Let us take $\psi $ as in (\ref{example 11}), then for any $k\ge 1$,
\begin{equation*}
\begin{array}{llll}
		\psi(2k)&=e^{-[\ln(2k)]^2}=e^{-(\ln2+\ln k)^2}=e^{-(\ln k)^2-2\ln k\ln2-(\ln2)^2}\\
		&=e^{-(\ln k)^2}e^{-2\ln k\ln2-(\ln2)^2}=\psi(k)e^{-(\ln2)(2\ln k+\ln2)}\\
		&=\psi(k)e^{-(\ln2)[\ln(2k^2)]}=\psi(k)e^{\ln(2k^2)^{-\ln2}}=\psi(k)(2k^2)^{-\ln2}\\
		&=\psi(k)\left( \frac{1}{2k^2}\right) ^{\ln2} \le \frac {1}{ 2^{\ln 2}} \frac {(2k)^A\psi (k) +\psi (k)} {k ^{2\ln 2}}.
\end{array}	
\end{equation*}
This shows that (\ref{the basic inequality 3}) holds true with $k_0=1$, $B=C=1$, $c_2=\frac 1 {2\ln 2}$, $D=2\ln 2$, and any $0< A< 2\ln2$.

Now we are going to show that (\ref{the basic inequality  2}) does not hold true with $B=C=1$, for any choice of the three constants $D>A>0$ and $c_1>0$: by contradiction, if (\ref{the basic inequality  2}) would hold true, then Lemma \ref{generalization of Stampacchia Lemma}  ii) would guarantees (\ref{B=C=1}), then (note that $\psi (k_0) =\psi (1)=1$)
\begin{equation*}
\psi(k)\leq e ^{1- \left(\frac {k-1}{\tau}\right) ^{\frac {D-A}{D}}},~~~ \forall k\in[1,+\infty)
\end{equation*}
for a suitable constant $\tau$ depending only on the constants $c_1, A$ and $D$. That is
\begin{equation*}
e^{-(\ln k)^2}\leq e ^{1- \left(\frac {k-1}{\tau}\right) ^{\frac {D-A}{D}}},
\end{equation*}
this means that
\begin{equation*}
 1 \leq e ^{1- \left(\frac {k-1}{\tau}\right) ^{\frac {D-A}{D}} +(\ln k)^2},
\end{equation*}
but this is false since
$$
\lim_{k\rightarrow +\infty }  1- \left(\frac {k-1}{\tau}\right) ^{\frac {D-A}{D}} +(\ln k)^2 =-\infty.
$$
\qed

\begin{rem} For the case $B>C>1$, one has
	\begin{equation*}\label{R2}
		(\ref{the basic inequality  2})\nLeftarrow(\ref{the basic inequality 3}).
	\end{equation*}
More precisely, the function
\begin{equation}\label{example 2}
	\psi(k)=e^{-k^{p}},~~~p=\log_{2}(2C),~~~k\in[1,+\infty)
\end{equation}
verifies (\ref{the basic inequality 3}) with $B>C>1$, $c_2=1$, any $D>0$ and a suitable $k_0=k_0(D,C)\geq 1$, but it does not satisfy (\ref{the basic inequality  2}) for any choice of the four constants $B>C>1$, $D>0$ and $c_1>0$.
\end{rem}
\noindent {\bf Proof. }
Let us take $\psi$ as in (\ref{example 2}), we keep in mind that $2^p=2C$ and we have
\begin{equation*}
\psi(2k)=e^{-(2k)^p}=e^{-2^p k^p}=e^{-2Ck^p}=(e^{-k^p})^{2C}=\left( \psi(k)\right) ^{2C}=( e^{-k^p}) ^C\left( \psi(k)\right) ^{C}.
\end{equation*}
Note that there exists $k_0=k_0(D,C)\geq 1$ such that
\begin{equation*}
( e^{-k^p}) ^C\leq \left( \frac{1}{k}\right) ^{D},~~~\forall k\in[k_0,+\infty).
\end{equation*}
Then
\begin{equation*}
\psi(2k)=( e^{-k^p}) ^C\left( \psi(k)\right) ^{C}\leq \left(  \frac{1}{k}\right) ^{D}\left( \psi(k)\right) ^{C},~~~\forall k\in[k_0,+\infty),
\end{equation*}
so that $\psi$ verifies (\ref{the basic inequality 3}) with $c_2=1$, with any $D>0$ and with a suitable $k_0=k_0(D,C)\geq 1$. We claim that such a $\psi$ does not satisfy (\ref{the basic inequality  2}) for any choice of the constants $B>C>1$, $D>0$, $c_1>0$, $k_0\geq 1$. Indeed, if such a $\psi$ would satisfy (\ref{the basic inequality  2}), then part iii) of  Lemma \ref{generalization of Stampacchia Lemma} would imply
(\ref{B>C>1}):
\begin{equation*}
\psi(2L)=0
\end{equation*}
for a suitable $L\geq0$: this gives a contradiction since $\psi(k)>0$ for every $k\in[1,+\infty)$.   \qed

\section{An Application.}
In this section, we shall answer the two questions proposed at the end of the first section.




\begin{theorem}\label{theorem}
	Let $f$ be in $M^{r}(\Omega)$, with $r\geq \frac n p$, and $u$ is a minima of $\cal I$ on $W_{0}^{1,q}(\Omega)$. Then
	
	(i) $r=\frac{n}{p}$ $\Rightarrow$ $\exists\lambda>0$ such that $e^{\lambda|u|^{1-p'\alpha}}\in L^{1}(\Omega)$;
	
	(ii) $r>\frac{n}{p}$ $\Rightarrow$ $\exists L>0$ such that $|u|\leq 2L$, a.e. $\Omega$.
	
\end{theorem}
\noindent {\bf Proof. } As in the proof of Theorem 6.3 in \cite{BO}, for $u\in W_0^{1,q} (\Omega)$ a minima of $\cal I$, we take
$$
v=T_k(u) =\min \{-k, \max \{k,u\}\}
$$
in (\ref{definition for minima}) and we get for any $k>0$,
$$
\int_{A_k} a(x,u) j(\nabla u) dx \le \int_{A_k} fG_k(u)dx,
$$
where
$$
A_k =\{x\in \Omega: |u(x)| \ge k\}, \ \ G_k(u) =u-T_k(u).
$$
We follow the lines of the proof of Theorem 6.3 in \cite{BO} until we arrive at the following inequality: for any $h>k>0$,
$$
\begin{array}{llll}
\displaystyle |A_h| &\leq &\displaystyle \frac{c}{(h-k)^{q^{*}}}\left[ |A_k|^{(p-1-\frac{q}{r}+\frac{q}{n})\frac{q^*}{q(p-1)}}k^{\frac{\alpha pq^*}{p-1}}+|A_k|^{(q-1-\frac{q}{r}+\frac{q}{n})\frac{q^*}{q[p(1-\alpha)-1]}}\right] \\
&\leq & \displaystyle \frac{c}{(h-k)^{q^{*}}}\left[ |A_k|^{(p-1-\frac{q}{r}+\frac{q}{n})\frac{q^*}{q(p-1)}}h^{\frac{\alpha pq^*}{p-1}}+|A_k|^{(q-1-\frac{q}{r}+\frac{q}{n})\frac{q^*}{q[p(1-\alpha)-1]}}\right],
\end{array}
$$
where $q$ is as in (\ref{definitioin for q}). We now apply Lemma \ref{generalization of Stampacchia Lemma} with
$$
\psi(k)=|A_k|, \ c_1=c, \  A=\frac{\alpha pq^*}{p-1},~~B=\left( p-1-\frac{q}{r}+\frac{q}{n} \right)\frac{q^*}{q(p-1)},
$$
$$
C=\left( q-1-\frac{q}{r}+\frac{q}{n} \right)\frac{q^*}{q[p(1-\alpha)-1]}, ~~D=q^* ~  \mbox { and }~ ~k_0=0.
$$
We note that $A<D$ since (\ref{condition for alpha No. 1}).

(i) If $r=\frac{n}{p}$, then $B=C=1$. We use Lemma \ref{generalization of Stampacchia Lemma} ii) and we derive that there exists a constant $\tau$ such that for any $k\ge 0$,
$$
|\{|u|\ge k\}|\leq |\{|u|\ge 0 \}|e^{1-\left( \frac{k}{\tau}\right)^{\frac{D-A}{D}}}\leq |\Omega| e^{1-\left( \frac{k}{\tau}\right)^{1-\alpha p'}} =|\Omega|e e ^{-2\lambda k^{1-\alpha p'}},
$$
where $2\lambda=\left( \frac{1}{\tau}\right) ^{1-\alpha p'}$ and we have used $\frac {D-A}{D}=1-\alpha p'$.

The above inequality implies
$$
|\{e^{\lambda|u|^{1-\alpha p' }}\ge e^{\lambda k^{1-\alpha p' }}\}|=|\{|u|\ge k\} \leq |\Omega|e e ^{-2\lambda k^{1-\alpha p'}}.
$$
Let $\tilde{k}=e^{\lambda k^{1-\alpha p' }}$, then
$$
|\{e^{\lambda|u|^{1-\alpha p' }}\ge \tilde{k}\}|\leq\frac{|\Omega|e}{\tilde{k}^{2}}, \ \  \forall \tilde{k}\geq 1.
$$
We now use Lemma 3.11 in \cite{Boccardo-Croce} which states that the sufficient and necessary condition for $g\in L^{r}(\Omega),r\geq1$, is
$$
\sum\limits_{k=1 }^{\infty}k^{r-1}|\{|g|\ge k\}|<+\infty.
$$
We use the above lemma for $g=e^{\lambda  |u|^{1-\alpha p' }}$ and $r=1$. Since
$$
\sum_{\tilde k =1} ^{\infty}  |\{ e^{\lambda  |u|^{1-\alpha p' }} \ge \tilde k  \}| \le |\Omega|e \sum_{\tilde k =1} ^{\infty} \frac 1 {\tilde k ^2}<+\infty,
$$
then $e^{\lambda  |u|^{1-\alpha p' }}\in L^{1}(\Omega)$, as desired.

(ii) If $r>\frac{N}{p}$, then $B>C>1$. We use Lemma \ref{generalization of Stampacchia Lemma} iii) and we have$ |A_{2L}|=0$ for some constant $L>0$, from which we derive $|u|\leq2L$ a.e. $\Omega$. That is, $u$ is bounded.
\qed

\vspace{3mm}

\noindent {\bf Acknowledgments:} The first author thanks NSFC(12071021), NSF of Hebei Province (A2019201120) and the Key Science and Technology Project of
Higher School of Hebei Province (ZD2021307) for the support; the second author thanks the Postgraduate Innovation
Project of Hebei Province (CXZZSS2020005).

\rm 
\end{document}